# On the sandpile group of regular trees

Evelin Toumpakari[*]

March 12, 2004

## Abstract

The sandpile group of a connected graph is the group of recurrent configurations in the abelian sandpile model on this graph. We study the structure of this group for the case of regular trees. A description of this group is the following: Let $\mathcal{T}(d,h)$ be the $d$-regular tree of depth $h$ and let $V$ be the set of its vertices. Denote the adjacency matrix of $\mathcal{T}(d,h)$ by $A$ and consider the modified Laplacian matrix $\Delta := dI - A$. Let the rows of $\Delta$ span the lattice $\Lambda$ in $\mathbb{Z}^V$. The sandpile group of $\mathcal{T}(d,h)$ is $\mathbb{Z}^V/\Lambda$. We compute the rank, the exponent and the order of this abelian group and find a cyclic Hall-subgroup of order $(d-1)^h$. We find that the base $(d-1)$-logarithm of the exponent and of the order are asymptotically $3h^2/\pi^2$ and $c_d(d-1)^h$, respectively. We conjecture an explicit formula for the ranks of all Sylow subgroups.

## 1 Introduction

Motivated by the concept of "self-organized criticality" in statistical mechanics [1] (cf. [10]), the Abelian Sandpile Model (ASM) is a game on a connected graph with a special vertex called the "sink." We refer to this graph as the "augmented graph" and the term "the graph" will be reserved to the augmented graph with the sink deleted. Vertices other than the sink will be called "ordinary." A configuration of the game is an assignment of a non-negative integer $w_i$ to each vertex $i$ of the augmented graph. These integers may be thought of as the numbers of sandgrains at the "sites." A configuration is stable if for all ordinary vertices $i$, $0 \leq w_i < \deg(i)$, where $\deg(i)$ is the degree of $i$ in the augmented graph. When the height at an ordinary vertex $i$ reaches $\deg(i)$, a "toppling" occurs, i.e., this vertex loses $\deg(i)$ grains, one to each of its neighbors (including the sink). Starting with any unstable configuration and toppling unstable ordinary vertices repeatedly, we finally arrive at a stable configuration. The order in which the topplings occur does not matter [4], therefore the model is called abelian. Moreover, the connectedness of the augmented graph ensures that a stable configuration will be reached in a finite number of steps: the sink "collects" the grains "falling off" the ordinary vertices. For every ordinary vertex $i$, we define an operator $\alpha_i$ on the space of stable configurations. For a stable configuration $\mathbf{w}$, $\alpha_i(\mathbf{w})$ is the stable configuration we obtain after adding one grain at vertex $i$ and toppling if necessary.

[*]Department of Mathematics, University of Chicago, 5734 S.University ave, Chicago, IL 60637, Email: evelint@math.uchicago.edu

A configuration $\mathbf{w}$ is *recurrent* if for every operator $\alpha_i$ there is a positive integer $r_i$ such that $\alpha_i^{r_i}(\mathbf{w}) = \mathbf{w}$ (Dhar [7]). The $\alpha_i$ restricted to the set of recurrent configurations form an abelian group of order equal to the number of spanning trees of the graph [7]. The operators $\alpha_i$ on the set of recurrent configurations satisfy the relations $\alpha_i^{\deg(i)} = \prod \alpha_j$, where the product extends over all ordinary neighbors $j$ of $i$ [7]. It turns out that these relations define the group.

Creutz [6] proved that the recurrent configurations under the operation of pointwise addition and toppling generate an abelian group isomorphic to the group generated by the operators $\alpha_i$. This group is called in [5] the *sandpile group* of the augmented graph.

Of special interest are sequences of sandpile groups derived from infinite $d$-regular graphs in the following way: Let $\mathcal{T} = (V(\mathcal{T}), E(\mathcal{T}))$ be an infinite $d$-regular graph. With every finite subset $P \subset V(\mathcal{T})$ we associate a graph $\hat{\mathcal{P}}$ as follows: We consider the subgraph $\mathcal{P}$ induced on $P$, adjoin a sink to $P = V(\mathcal{P})$ and join each vertex $i \in P$ by $d - \deg_{\mathcal{P}}(i)$ edges to the sink. This way all vertices in $P$ have degree $d$ in $\hat{\mathcal{P}}$. Now take a "well-behaved" sequence $\{P_n\}$ of subsets of $V(\mathcal{T})$ and consider the sandpile groups $G_n$ associated with $\hat{\mathcal{P}}_n$. An instance of this procedure is the sequence of finite square lattices studied in [9]. With some abuse of terminology, we shall refer to $G$ as the sandpile group of the graph $\mathcal{P}$ rather than the sandpile group of the augmented graph $\hat{\mathcal{P}}$.

In this paper we choose $\mathcal{T}$ to be the infinite regular tree of degree $d$ and $P_h$ to be the ball of radius $h$ about a vertex designated as the root. We study the associated sandpile group $G(d, h)$, which we call the *sandpile group of the $d$-valent tree of depth $h$*. We find the rank, the exponent and the order of $G(d, h)$ and we show that $G(d, h)$ has a cyclic Hall-subgroup of order $(d-1)^h$.

A slight variation of the ASM is studied in [2] and [3]. Probabilistic aspects of the ASM on $\hat{\mathcal{T}}(d, h)$ and on infinite regular trees are studied in [8] and [11] respectively. A comprehensive introduction to the ASM and to similar models can be found in [10].

**Acknowledgments:** I would like to thank my advisor Steven Lalley for introducing me to this topic and for his constant encouragement. I would like to express my gratitude to László Babai for many helpful discussions and for his support.

## 2  Main results

In this section we state the main results of this paper. The results concern $G(d, h)$, the sandpile group of the $d$-valent tree of depth $h$. Throughout this paper we assume $d \geq 3$, $h \geq 1$.

**Theorem 2.1** *The* rank *of $G(d,h)$ is $(d-1)^h$.*

**Definition 2.2** Let $G$ be a finite group. The *exponent* of $G$ is the least common multiple of the orders of the elements of $G$. If $G$ is abelian then the exponent is also the largest order of an element.

**Notation:**

(i) We denote the exponent of $G(d,h)$ by $\exp(d,h)$.

(ii) We define the numbers $\theta(d,n)$ as $\theta(d,n) := [(d-1)^n - 1]/(d-2)$.

**Theorem 2.3** *The exponent $\exp(d,h)$ of the group $G(d,h)$ is equal to*

$$(d-1)^h \mathrm{lcm}\{d\theta(d,h+1), \theta(d,h), \theta(d,h-1), \ldots, \theta(d,2)\}. \tag{1}$$

This seems to be the first time in the literature that substantial information about the exponent of a nontrivial class of sandpile groups is found.

**Corollary 2.4** *For every fixed $d \geq 3$, the following asymptotic equality holds as $h \to \infty$:*

$$\log_{d-1} \exp(d,h) \sim \frac{3h^2}{\pi^2}.$$

**Notation:** We denote the order of $G(d,h)$ by $g(d,h)$.

**Theorem 2.5** *The order $g(d,h)$ of the group $G(d,h)$ is equal to*

$$d(d-1)^h [\theta(d,h+1)]^{d-1} \prod_{n=1}^{h-1} [\theta(d,h+1-n)]^{(d-2)d(d-1)^{n-1}}. \tag{2}$$

**Corollary 2.6** *For every fixed $d \geq 3$, the following asymptotic equality holds for $g(d,h)$ as $h \to \infty$:*

$$\log_{d-1} g(d,h) \sim c_d (d-1)^h, \tag{3}$$

*where*

$$c_d := d(d-2) \sum_{n=0}^{\infty} (d-1)^{-2-n} \log_{d-1} \frac{(d-1)^{n+2} - 1}{d-2}.$$

**Definition 2.7**

(i) Let $G$ be a group and $H$ a subgroup of $G$. $H$ is a *Hall-subgroup* if $\gcd(|H|, |G:H|) = 1$.

(ii) We say that $H$ is a *Hall $t$-subgroup* of $G$ if $H$ is a Hall-subgroup of $G$ and $t$ and $|H|$ have the same prime divisors. E.g., Hall 6-subgroups are the same as Hall 18-subgroups.

**Note:** If $G$ is abelian then for any $t \mid |G|$, $G$ has a unique Hall $t$-subgroup.

**Theorem 2.8** *The Hall $(d-1)$-subgroup of $G(d,h)$ is cyclic of order $(d-1)^h$.*

# 3 Sandpile groups: general theory

Let $\hat{\mathcal{R}}$ be a finite graph with a vertex designated as the sink and let $\mathcal{R} = (V(\mathcal{R}), E(\mathcal{R}))$ be the subgraph of $\hat{\mathcal{R}}$ obtained after deletion of the sink. Let $A$ denote the adjacency matrix of $\hat{\mathcal{R}}$ and $D$ denote the diagonal matrix with $D_{ii} := \deg_{\hat{\mathcal{R}}}(i)$. The *Laplacian* of $\hat{\mathcal{R}}$ is defined to be the matrix $L := D - A$. Let $\Delta$ be the submatrix of $L$ obtained after deletion of the row and column corresponding to the sink and let $\boldsymbol{\delta}_i$ ($i \in V(\mathcal{R})$) be the row vectors of $\Delta$. Let $\{\mathbf{x}_i : i \in V(\mathcal{R})\}$ be the standard basis in $\mathbb{Z}^{V(\mathcal{R})}$. Also, for $i \in V(\mathcal{R})$ let $N_i$ denote the set of neighbors of $i$ in $\mathcal{R}$. By the definition of $\Delta$ it is clear that

$$\boldsymbol{\delta}_i = \deg_{\hat{\mathcal{R}}}(i)\mathbf{x}_i - \sum_{j \in N_i} \mathbf{x}_j \quad (i \in V(\mathcal{R})). \tag{4}$$

The relations in the group of operators $\alpha_i$ mentioned in the introduction can be written in terms of the Laplacian as

$$\prod_{j \in V(\mathcal{R})} \alpha_i^{\Delta_{ij}} = 1 \ (i \in V(\mathcal{R})). \tag{5}$$

Consequently, the sandpile group $G$ of $\hat{\mathcal{R}}$ is expressed as

$$G := \mathbb{Z}^{V(\mathcal{R})}/\Lambda, \tag{6}$$

where $\Lambda := \sum_{i \in V(\mathcal{R})} \mathbb{Z}\boldsymbol{\delta}_i$ is the lattice spanned by the $\boldsymbol{\delta}_i$.

**Notation:** Let $G$ be a group and let $\mathbf{z} \in G$. We denote the order of $\mathbf{z}$ in $G$ by $\mathrm{ord}(\mathbf{z})$.

The following lemma will be used repeatedly:

**Lemma 3.1** *Let $\mathbf{u}_1, \ldots, \mathbf{u}_t \in \mathbb{Z}^t$ be linearly independent over $\mathbb{Q}$ and let $\Lambda := \sum_{n=1}^t \mathbb{Z}\mathbf{u}_n$ be the lattice spanned by the $\mathbf{u}_n$. Consider the finite abelian group $K = \mathbb{Z}^t/\Lambda$. Assume $\mathbf{u} \in \mathbb{Z}^t$ satisfies*

$$r\mathbf{u} = \sum_{n=1}^t r_n \mathbf{u}_n, \tag{7}$$

*where $r, r_n \in \mathbb{Z}$, $r > 0$ and $\gcd(r_1, \ldots, r_n) = 1$. Then the order of $\mathbf{u}$ in $K$ is $r$.*

**Proof:** Let $\mathbf{u} \mapsto \bar{\mathbf{u}}$ denote the quotient map $\mathbb{Z}^t \to \mathbb{Z}^t/\Lambda$. By definition we have $r\mathbf{u} \in \Lambda$, therefore $\mathrm{ord}(\bar{\mathbf{u}}) \mid r$. Let $r = \mathrm{ord}(\bar{\mathbf{u}})s$. The relation $\mathrm{ord}(\bar{\mathbf{u}})\mathbf{u} \in \Lambda$ implies

$$\mathrm{ord}(\bar{\mathbf{u}})\mathbf{u} = \sum_{n=1}^t q_n \mathbf{u}_n, \tag{8}$$

where $q_n \in \mathbb{Z}$. Myltiplying (8) by $s$ we obtain

$$r\mathbf{u} = \sum_{n=1}^{t} sq_n \mathbf{u}_n. \qquad (9)$$

The $\mathbf{u}_n$ are linearly independent and thus equations (7) and (9) yield $r_n = sq_n$ ($1 \leq n \leq t$). So $s \mid \gcd(r_1, \ldots, r_n)$ and therefore $s = 1$. □

**Definition 3.2** Let $s$ be a positive integer and let $p$ be a prime.

(i) We define $e_p(s) := w$ if $p^w \mid s$ and $p^{w+1} \nmid s$.

(ii) Let $G$ be group and $\mathbf{z} \in G$. We define $e_p(\mathbf{z}) := e_p(\text{ord}(\mathbf{z}))$.

The following Lemma will be used for proving Theorem 2.3.

**Lemma 3.3** Let $G$ be an abelian group and let $\mathbf{z}_n \in G$ ($0 \leq n \leq t$). Let $\text{ord}(\mathbf{z}_0) = s$. Assume $r_n \mathbf{z}_{n+1} = r_{n+1} \mathbf{z}_n$ ($0 \leq n \leq t-1$), where $r_n \in \mathbb{Z}$ and $\gcd(r_n, r_{n+1}) = 1$ ($0 \leq n \leq t-1$). Then $\text{ord}(\mathbf{z}_n) \mid \text{lcm}\{sr_0, r_1, \ldots, r_{n-1}\}$ ($1 \leq n \leq t$).

**Proof:** Let $p$ be a prime.

**Claim 3.4**

$$e_p(\mathbf{z}_n) \leq \max\{e_p(\mathbf{z}_{n-2}), e_p(\mathbf{z}_{n-1}), e_p(r_{n-1})\} \ (2 \leq n \leq t). \qquad (10)$$

**Proof of Claim 3.4:** We have

$$r_{n-1}\mathbf{z}_n = r_n \mathbf{z}_{n-1}. \qquad (11)$$

If $p \nmid r_{n-1}$ then by equation (11) we have $e_p(\mathbf{z}_n) \leq e_p(\mathbf{z}_{n-1})$.
If $p \mid r_{n-1}$ consider the equation

$$r_{n-2}\mathbf{z}_{n-1} = r_{n-1}\mathbf{z}_{n-2}. \qquad (12)$$

We have $p \nmid r_{n-2} r_n$.
If $e_p(r_{n-1}) < e_p(\mathbf{z}_{n-2})$, then by equation (12) we obtain

$$e_p(\mathbf{z}_{n-1}) = e_p(\mathbf{z}_{n-2}) - e_p(r_{n-1}) > 0. \qquad (13)$$

So by equation (11) we obtain $e_p(\mathbf{z}_n) = e_p(\mathbf{z}_{n-1}) + e_p(r_{n-1})$. By (13) this is $e_p(\mathbf{z}_{n-2})$.

If $e_p(r_{n-1}) \geq e_p(\mathbf{z}_{n-2})$ then by equation (12) we obtain $e_p(\mathbf{z}_{n-1}) = 0$ and by equation (11) we obtain $e_p(\mathbf{z}_n) \leq e_p(r_{n-1})$. □

For the proof of the Lemma, it suffices to prove the following:

$$e_p(\mathbf{z}_n) \leq \max\{e_p(s) + e_p(r_0), e_p(r_1), \ldots, e_p(r_{n-1})\} \ (1 \leq n \leq t). \qquad (14)$$

We prove (14) by induction on $n$:

Let $n = 1$. We have
$$r_0 \mathbf{z}_1 = r_1 \mathbf{z}_0 \Rightarrow r_0 s \mathbf{z}_1 = r_1 s \mathbf{z}_0 = 0.$$

Let $n = 2$. By the Claim we have
$$e_p(\mathbf{z}_2) \leq \max\{e_p(\mathbf{z}_0), e_p(\mathbf{z}_1), e_p(r_1)\} \leq \max\{e_p(s) + e_p(r_0), e_p(r_1)\}.$$

For $n \geq 3$, the inductive step is accomplished by Claim 3.4. $\square$

## 4 Preliminaries and notation

From now on $\mathcal{T} = \mathcal{T}(d) = (V(\mathcal{T}(d)), E(\mathcal{T}(d)))$ denotes the infinite $d$-valent tree. The distance of the vertices $i$ and $j$ is denoted by $\mathrm{dist}(i,j)$. We choose a vertex to be the root and we designate it by 0. Let $i, j \in V(\mathcal{T}(d))$. We say that $i$ is an *ancestor* of $j$ or $j$ is a *descendant* of $i$ ($i \preceq j$), if $i$ is on the unique path from $j$ to 0 and that $j$ is a *proper descendant* of $i$ ($i \prec j$), if $i \preceq j$ and $i \neq j$. Let $V_i = \{k : i \preceq k\}$ be the set of descendants of $i$. We say that $i$ is the *parent* of $j$ ($i = p(j)$) or that $j$ is a *child* of $i$, if $i \prec j$ and $\mathrm{dist}(i,j) = 1$. Let $C_i = \{k : p(k) = i\}$ be the set of children of $i$. In particular, let $C_0 = \{1, \ldots, d\}$ denote the set of children of the root 0. We say that $i$ and $j$ are *siblings* if $p(i) = p(j)$ and $i \neq j$. Let $B_n$ be the ball of radius $n$ about the root, i.e. $B_n := \{k \in V(\mathcal{T}(d)) : \mathrm{dist}(0,k) \leq n\}$ and let $S_n$ be the sphere of radius $n$ about the root, i.e. $S_n := \{k \in V(\mathcal{T}(d)) : \mathrm{dist}(0,k) = n\}$. Let $V_i^n := V_i \cap S_n$ be the set of descendants of $i$ in depth $n$.

We define the $d$-valent tree of depth $h$ $\mathcal{T}(d,h) = (V(\mathcal{T}(d,h)), E(\mathcal{T}(d,h)))$ as the subgraph of $\mathcal{T}(d)$ induced on $B_h$. A vertex $i$ is called a *leaf* if $\mathrm{dist}(0,i) = h$. In our notation $V_i^h = V_i \cap S_h$ is the set of leaves below the vertex $i$. Let $\hat{\mathcal{T}}(d,h)$ be the graph obtained from $\mathcal{T}(d,h)$ by attaching an additional vertex (the sink) by $d-1$ edges to each leaf. We denote the sandpile group of $\hat{\mathcal{T}}(d,h)$ by $G(d,h)$ and call it the *sandpile group of the $d$-valent tree of depth $h$*. The order of an element $\mathbf{x}$ in $G(d,h)$ is denoted by $\mathrm{ord}(\mathbf{x})$.

In the case of $G(d,h)$ expression (4) turns into
$$\boldsymbol{\delta}_i := \begin{cases} d\mathbf{x}_i - \mathbf{x}_{p(i)} - \sum_{j \in C_i} \mathbf{x}_j & \text{if } \mathrm{dist}(0,i) < h \\ d\mathbf{x}_i - \mathbf{x}_{p(i)} & \text{otherwise} \end{cases}. \quad (15)$$

The numbers $\theta(d,n)$ defined as follows occur many times throughout the analysis:
$$\theta(d,n) := \frac{(d-1)^n - 1}{d-2}. \quad (16)$$

The $\theta(d,n)$ satisfy the recurrence
$$\theta(d, n+2) = d\theta(d, n+1) - (d-1)\theta(d, n). \quad (17)$$

# 5 Table of notation

**Graphs**

$\mathcal{T}(d) = (V(\mathcal{T}(d)), E(\mathcal{T}(d)))$: infinite $d$-regular tree.
$\mathcal{P}$: subgraph induced on $P \subset V(\mathcal{T}(d))$.
$\hat{\mathcal{P}}$: $\mathcal{P}$ with sink adjoined, such that all vertices in $P$ have degree $d$ in $\hat{\mathcal{P}}$.
$\mathcal{T}(d,h)$: $d$-regular tree of depth $h$.
$\hat{\mathcal{T}}(d,h)$: $\mathcal{T}(d,h)$ with sink adjoined.
$\hat{\mathcal{R}}$: finite graph with a vertex designated as the sink.
$\mathcal{R} = (V(\mathcal{R}), E(\mathcal{R}))$: subgraph of $\hat{\mathcal{R}}$ obtained after deletion of the sink.
$\mathcal{S}_j(d,h)$: the subgraph of $\mathcal{T}(d,h)$ induced on $V_j$ (see Definition 8.2).
$\hat{\mathcal{S}}_j(d,h)$: $\mathcal{S}_j(d,h)$ with "sink" adjoined (see Definition 8.2).

**Vertices**

$i, j, k, m$
$0$: the root
$1, \ldots, d$: the vertices at depth 1.
$m_i$: some child of the vertex $i$.
$\ell$: a leaf
$\ell_i$: some leaf that is descendant of the vertex $i$.

**Sets of vertices**

$P$: finite subset of $V(\mathcal{T}(d))$.
$P_n$: sequence of subsets of $V(\mathcal{T}(d))$.
$N_i$: set of neighbors of the vertex $i$.
$V_i$: set of descendants of $i$.
$U_i := V_i \backslash \{\ell_i\}$.
$C_i$: set of children of $i$.
$J_i := C_i \backslash \{m_i\}$.
$B_n$: ball of radius $n$ about the root.
$S_n$: sphere of radius $n$ about the root.
$V_i^n := V_i \cap S^n$: set of vertices at depth $n$ that are descendants of the vertex $i$.
$V_i^h$: set of leaves of $\mathcal{T}(d,h)$ that are descendants of the vertex $i$.
$V := V(\mathcal{T}(d,h))$: the set of vertices of $\mathcal{T}(d,h)$.
$F(d,h)$: subset of $V(\mathcal{T}(d,h))$ (see Definition 6.2).
$B$: subset of $V(\mathcal{T}(d,h))$.

**Constants**

$d$: degree of the regular tree $\mathcal{T}(d)$.
$h$: depth of finite tree $\mathcal{T}(d,h)$.
$\exp(d,h)$: the exponent of $G(d,h)$.

$f(d, h)$: the order of $\mathbf{x}_\ell$ for $\ell$ leaf.
$g(d, h)$: the order of $G(d, h)$.
$\theta(d, n) := [(d-1)^n - 1]/(d-2)$.
$\eta(r, s) := \operatorname{lcm}\{r^n - 1 : n = 1, \ldots, s\}$.
$c_d$: a constant (see Corollary 2.6).

**Integers**

$p, q, r, s, t$.
$n$: integer usually between 0 and $h$.

**Matrices**

$A$: adjacency matrix of graph $\hat{\mathcal{R}}$ with sink.
$D$: diagonal matrix with $D_{ii} := \deg_{\hat{\mathcal{R}}}(i)$.
$L := D - A$: Laplacian of graph $\hat{\mathcal{R}}$ with sink.
$\Delta$: matrix obtained from $L$ after deleting row and column corresponding to the sink.

**Groups**

$G_n$: sandpile group of $\hat{\mathcal{P}}_n$.
$G(d, h)$: sandpile group associated with $B_h$.
$G$: sandpile group of $\hat{\mathcal{R}}$.
$\Lambda$: lattice in $\mathbb{Z}^t$ spanned by $t$ linearly independent elements of $\mathbb{Z}^t$.
$K := \mathbb{Z}^t/\Lambda$.
$G_F(d, h)$: subgroup of $G(d, h)$ generated by $\{\bar{\mathbf{x}}_i : i \in F(d, h)\}$.
$G_n(d, h)$: subgroup of $G(d, h)$ generated by $\{\bar{\mathbf{x}}_i : i \in B_n\}$.
$M_j(d, h)$: sandpile group of $\hat{\mathcal{S}}_j(d, h)$ (see Definition 8.2).
$K_j(d, h) := M_j(d, h)/\langle \bar{\mathbf{x}}_j \rangle$ (see Definition 8.2).
$H$: Hall subgroup of a group $G$.
$S_p(d, h)$: the Sylow $p$-subgroup of $G(d, h)$.

**Elements of groups/vectors**

$\boldsymbol{\delta}_i$: the row vector of $\Delta$ corresponding to the vertex $i$.
$\mathbf{x}_i$: the standard basis vector of $\mathbb{Z}^{V(\mathcal{T}(d,h))}$ corresponding to the vertex $i$.
$\mathbf{v} = (v_i)_{i \in V(\mathcal{T}(d,h))}$: an element of $\mathbb{Z}^{V(\mathcal{T}(d,h))}$.
$\bar{\mathbf{v}} := \phi(\mathbf{v})$.
$\mathbf{u}_1, \ldots, \mathbf{u}_t$: linearly independent elements of $\mathbb{Z}^t$.
$\mathbf{u}$: some element in the lattice $\sum_{n=1}^{t} \mathbb{Z}\mathbf{u}_n$ generated by the $\mathbf{u}_n$.
$\bar{\mathbf{u}}$: the image of $u$ in the quotient group $\mathbb{Z}^t/\sum_{n=1}^{t} \mathbb{Z}\mathbf{u}_n$.
$\mathbf{s}_h = (s_{h,i})_{i \in V(\mathcal{T}(d,h))}$: a solution of system (19).
$\bar{\mathbf{y}}_n = \bar{\mathbf{y}}(d, n) := \sum_{i \in S_n} \bar{\mathbf{x}}_i$: element of $G(d, h)$.

$\mathbf{z}, \mathbf{z}_t$: elements of an abelian group $G$.

**Vector spaces**

$\Gamma(j, B)$: the vector subspace of $\mathbb{Q}^B$ generated by the $\{\mathrm{pr}_B(\boldsymbol{\delta}_k) : k \in V_j\}$.

**Functions**

$w_i$: height of sandpile at vertex $i$.
$\mathbf{w} = (w_i)_{i \in V(\mathcal{T}(d,h))}$: height configuration.
$\alpha_i$: toppling operator corresponding to the vertex $i$.
$p(i)$: the parent of the vertex $i$.
$\omega : \mathbb{Z}^{V(\mathcal{T}(d,h))} \to \mathbb{Z}_d^{F(d,h)}$.
$\phi$: quotient map from $\mathbb{Z}^{V(\mathcal{T}(d,h))}$ to $\mathbb{Z}^{V(\mathcal{T}(d,h))}/\sum_{i \in V(\mathcal{T}(d,h))} \mathbb{Z}\boldsymbol{\delta}_i$.
$\deg_{\mathcal{P}}(i)$: degree of the vertex $i$ in $\mathcal{P}$.
$\mathrm{dist}(i, j)$: distance of the vertices $i, j$ in $\mathcal{T}(d)$.
$\mathrm{ord}(\mathbf{z})$: the order of $\mathbf{z} \in G$ in the group $G$.
$e_p(s) := w$ if $p^w \mid s$ and $p^{w+1} \nmid s$.
$e_p(\mathbf{z}) := e_p(\mathrm{ord}(\mathbf{z}))$.
$\mathrm{pr}_B(\mathbf{v}) := (v_i)_{i \in B}$.
$\tilde{\mathbf{v}} := (v_i)_{i \in F(d,h)}$, i.e., $\tilde{\mathbf{v}} = \mathrm{pr}_{F(d,h)}(\mathbf{v})$.
$\mathrm{ord}_s(r)$: the multiplicative order of $r \mod s$ ($\gcd(r, s) = 1$ and $s > 0$).
$t_p(d)$: the least positive integer $n$ such that $p \mid \theta(d, n)$ ($p \nmid d - 1$).

# 6 The rank of the sandpile group $G(d, h)$

In this section we prove that the rank of $G(d, h)$ is $(d-1)^h$ (Theorem 2.1).

First, we introduce some notation and define a subset $F(d, h)$ of the vertex set of $\mathcal{T}(d, h)$, with $|F(d, h)| = (d - 1)^h$. We shall see that the generators $\{\bar{\mathbf{x}}_i : i \in F(d, h)\}$ generate $G(d, h)$ (Lemma 6.4). Recall that we denote the root by 0 and that $C_i$ is the set of children of the vertex $i$.

**Definition 6.1** Let $i \in V(\mathcal{T}(d, h))$. Let us pick a child $m_i$ of $i$. Set $J_i := C_i \backslash \{m_i\}$.

**Definition 6.2** If $h$ is even, then let

$$F(d, h) := \{0\} \bigcup \left( \bigcup_{q=0}^{(h-2)/2} \bigcup_{i \in S_{2q+1}} J_i \right).$$

If $h$ is odd, then let

$$F(d, h) := \bigcup_{q=0}^{(h-1)/2} \bigcup_{i \in S_{2q}} J_i.$$

**Notes:**

(i) If $h$ even, then $|F(d,h)| = 1 + (d-2) \sum_{q=0}^{(h-2)/2} d(d-1)^{2q} =$
$= 1 + d(d-2)[(d-1)^h - 1]/[(d-1)^2 - 1] = (d-1)^h$.

If $h$ odd, then $|F(d,h)| = (d-1) + (d-2) \sum_{q=1}^{(h-1)/2} d(d-1)^{2q-1} =$
$= (d-1) + (d-2)d(d-1) \sum_{q=0}^{(h-3)/2} (d-1)^{2q} =$
$= (d-1) + (d-2)d(d-1)[(d-1)^{h-1} - 1]/[(d-1)^2 - 1] =$
$= (d-1) + (d-1)[(d-1)^{h-1} - 1] = (d-1)^h$.

(ii) Observe that $F(d, h+2) = F(d, h) \cup (\cup_{i \in S_{h+1}} J_i)$.

**Notation:**

(i) Let $\Lambda := \sum_{i \in V(\mathcal{T}(d,h))} \mathbb{Z} \boldsymbol{\delta}_i$ denote the lattice generated by the $\delta_i$.

(ii) Let $\phi$ denote the quotient map from $\mathbb{Z}^{V(\mathcal{T}(d,h))}$ to $G(d,h) := \mathbb{Z}^{V(\mathcal{T}(d,h))}/\Lambda$.

(iii) Let $\mathbf{v} \in \mathbb{Z}^{V(\mathcal{T}(d,h))}$. We denote $\phi(\mathbf{v})$ by $\bar{\mathbf{v}}$.

**Definition 6.3** Let $G_F(d,h)$ be the subgroup of $G(d,h)$ generated by $\{\bar{\mathbf{x}}_i : i \in F(d,h)\}$.

The proof of Theorem 2.1 is based on the next two lemmas:

**Lemma 6.4** $G_F(d,h) = G(d,h)$.

**Proof:** In this proof we will work in $G(d,h)/G_F(d,h)$, so, e.g., the statement "$\bar{\mathbf{x}}_i = \bar{\mathbf{x}}_j$" means "$\bar{\mathbf{x}}_i + G_F(d,h) = \bar{\mathbf{x}}_j + G_F(d,h)$." Note that under this convention, $\bar{\mathbf{x}}_i = 0$ for all $i \in F(d,h)$. Also, note that by the definition of $G(d,h)$, $\bar{\boldsymbol{\delta}}_i = 0$ for all $i \in V(\mathcal{T}(d,h))$.

We prove by induction on $h$ that $G(d,h)/G_F(d,h) = 0$. We will show this by showing that $(\forall i \in V(\mathcal{T}(d,h)))(\bar{\mathbf{x}}_i = 0)$.

If $h = 0$, this is obvious.

If $h = 1$, then by the definition of $F(d,1)$ we have $(\forall i \in J_0)(\bar{\mathbf{x}}_i = 0)$. Also, $(\forall i \in J_0)(\bar{\boldsymbol{\delta}}_i = 0)$, so $(\forall i \in J_0)(d\bar{\mathbf{x}}_i - \bar{\mathbf{x}}_0 = 0)$, and therefore $\bar{\mathbf{x}}_0 = 0$. Finally, $\bar{\boldsymbol{\delta}}_0 = 0$, so $d\bar{\mathbf{x}}_0 - \sum_{i \in J_0} \bar{\mathbf{x}}_i - \bar{\mathbf{x}}_{m_0} = 0$, and thus $\bar{\mathbf{x}}_{m_0} = 0$.

Assume that the statement is true for $h$ and consider the case $h + 2$: We have $\cup_{i \in S_{h+1}} J_i \subset F(d, h+2)$, so $(\forall i \in S_{h+1})(\forall k \in J_i)(\bar{\mathbf{x}}_k = 0)$. Also, $(\forall i \in S_{h+1})(\forall k \in J_i)(\bar{\boldsymbol{\delta}}_k = 0)$, so $(\forall i \in S_{h+1})(\forall k \in J_i)(d\bar{\mathbf{x}}_k - \bar{\mathbf{x}}_i = 0)$, and hence $(\forall i \in S_{h+1})(\bar{\mathbf{x}}_i = 0)$. So, all relations of the form $\bar{\boldsymbol{\delta}}_i = 0$ with $i \in B_h$ for $G(d, h+2)/G_F(d, h+2)$ become identical with the corresponding ones for $G(d,h)/G_F(d,h)$, and therefore by the inductive assumption we obtain $(\forall i \in B_h)(\bar{\mathbf{x}}_i = 0)$. Finally, for $i \in S_{h+1}$ we have $\bar{\boldsymbol{\delta}}_i = 0$, therefore $d\bar{\mathbf{x}}_i - \bar{\mathbf{x}}_{p(i)} - \sum_{k \in J_i} \bar{\mathbf{x}}_k - \bar{\mathbf{x}}_{m_i} = 0$, and thus $\bar{\mathbf{x}}_{m_i} = 0$. □

**Corollary 6.5** The rank of the sandpile group $G(d,h)$ is $\leq (d-1)^h$.

**Lemma 6.6** $G(d,h)$ contains a copy of $\underbrace{\mathbb{Z}_d \oplus \cdots \oplus \mathbb{Z}_d}_{(d-1)^h \text{ times}}$.

**Proof:** We work in $\pmod{d}$ arithmetic. For $i \in F(d,h)$, fix $r_i \in \mathbb{Z}$ and consider the system

$$\begin{cases} \Delta \mathbf{v}^T = 0, & \mathbf{v} = (v_i)_{i \in V(\mathcal{T}(d,h))}, \\ v_i = r_i, & i \in F(d,h). \end{cases} \quad (18)$$

System (18) is equivalent to

$$\begin{cases} \boldsymbol{\delta}_i \cdot \mathbf{v} = 0, & i \in V(\mathcal{T}(d,h)), \\ v_i = r_i, & i \in F(d,h), \end{cases} \quad (19)$$

where $\mathbf{u} \cdot \mathbf{v}$ denotes the dot product of the vectors $\mathbf{u}$ and $\mathbf{v}$ $\pmod{d}$.

**Claim 6.7** For any choice of $r_i \in \mathbb{Z}$, system (19) has a solution in $\mathbb{Z}^{V(\mathcal{T}(d,h))}$.

**Proof of Claim 6.7:** By induction on $h$.

If $h = 0$, this is clear.

If $h = 1$, we have $F(d,1) = J_0$. Set $s_{1,0} := 0$, $s_{1,m_0} := -\sum_{i \in J_0} r_i$ and $s_{1,i} := r_i$ for $i \in J_0$. Then $\mathbf{s}_1 = (s_{1,i})_{i \in V(\mathcal{T}(d,1))}$ is a solution of (19).

Consider (19) for the tree of depth $h+2$ with chosen constants $r_i$ for $i \in F(d, h+2)$. Let $\mathbf{s}_h = (s_{h,i})_{i \in V(\mathcal{T}(d,h))}$ be a solution of (19) for the $h$ case with constants $r_i$ for $i \in F(d,h)$. Set $s_{h+2,i} := s_{h,i}$ for all $i \in B_h$, $s_{h+2,i} := 0$ for $i \in S_{h+1}$, $s_{h+2,m_i} = -\sum_{k \in J_i} r_k - s_{h,p(i)}$ for all $i \in S_{h+1}$ and $s_{h+2,i} := r_i$ for all $i \in F(d, h+2)$. Then $\mathbf{s}_{h+2} := (s_{h+2,i})_{i \in V(\mathcal{T}(d,h+2))}$ is a solution of (19). □

Claim 6.7 yields that for all $k \in F(d,h)$ there exists a solution vector $\mathbf{v}_k = (v_{k,i})_{i \in V(\mathcal{T}(d,h))}$ of (19) such that for all $k \in F(d,h)$ we have $v_{k,k} = 1$, and for all $k, i \in F(d,h)$ with $k \neq i$ we have $v_{k,i} = 0$.

For $\mathbf{v} \in \mathbb{Z}^{V(\mathcal{T}(d,h))}$, let $\tilde{\mathbf{v}} := (v_i)_{i \in F(d,h)}$. Note that $\{\tilde{x}_i : i \in F(d,h)\}$ is a basis of $\mathbb{Z}_d^{F(d,h)}$.

Consider the homomorphism $\omega : \mathbb{Z}^{V(\mathcal{T}(d,h))} \to \mathbb{Z}_d^{F(d,h)}$ with

$$\omega(\mathbf{u}) := (\mathbf{u} \cdot \mathbf{v}_k)_{k \in F(d,h)}.$$

For all $k \in F(d,h)$ we have $\omega(\mathbf{x}_k) = \tilde{\mathbf{x}}_k$. Therefore $\omega$ is surjective. Also, for all $k \in F(d,h)$, $\mathbf{v}_k$ is a solution of (19). Therefore $(\forall k \in F(d,h))(\forall i \in V(\mathcal{T}(d,h)))(\boldsymbol{\delta}_i \cdot \mathbf{v}_k = 0)$. So, $(\forall i \in V(\mathcal{T}(d,h)))(\omega(\boldsymbol{\delta}_i) = 0)$. Therefore there exists a homomorphism $\psi : G \to \mathbb{Z}_d^{F(d,h)}$ with $\omega = \psi \circ \phi$. The map $\omega$ is surjective therefore $\psi$ is surjective.

□

# 7 The exponent of $G(d, h)$

In this section we show that the exponent $\exp(d, h)$ of $G(d, h)$ is equal to $(d-1)^h \cdot \text{lcm} \{d\theta(d, h+1), \theta(d, h), \theta(d, h-1), \ldots, \theta(d, 2)\}$ (Theorem 2.3).

Recall that $\Lambda := \sum_{i \in V(\mathcal{T}(d,h))} \mathbb{Z} \boldsymbol{\delta}_i$ is the lattice spanned by the $\boldsymbol{\delta}_i$, that $\phi$ denotes the quotient map from $\mathbb{Z}^{V(\mathcal{T}(d,h))}$ to $G(d, h) := \mathbb{Z}^{V(\mathcal{T}(d,h))}/\Lambda$ and that $\bar{\mathbf{v}} := \phi(\mathbf{v})$. Also, recall that $V_j^q$ denotes the set of descendants of the vertex $j$, which are at depth $q$.

**Lemma 7.1** *Let $j \in S_n$ ($1 \leq n \leq h$).*
*Then $\theta(d, h+2-n)\bar{\mathbf{x}}_j - \theta(d, h+1-n)\bar{\mathbf{x}}_{p(j)} = 0$.*

**Proof:** Consider the sum

$$\sum_{q=n}^{h} [(d-1)^{h+1-q} - 1] \sum_{k \in V_j^q} \boldsymbol{\delta}_k. \tag{20}$$

Let $k \in V_j^q$ with $q > n$. The coefficient of $\mathbf{x}_k$ in sum (20) is

$$-(d-1)[(d-1)^{h-q} - 1] + d[(d-1)^{h+1-q} - 1] - [(d-1)^{h+2-q} - 1] =$$
$$-(d-1)^{h+1-q} + d - 1 + d(d-1)^{h+1-q} - d - (d-1)^{h+2-q} + 1 =$$
$$(d-1)(d-1)^{h+1-q} - (d-1)^{h+2-q} = 0.$$

The coefficient of $\mathbf{x}_j$ is

$$-(d-1)[(d-1)^{h-n} - 1] + d[(d-1)^{h+1-n} - 1] = (d-1)^{h+2-n} - 1.$$

Finally the coefficient of $\mathbf{x}_{p(j)}$ is

$$-[(d-1)^{h+1-n} - 1].$$

Therefore we have

$$\sum_{q=n}^{h} [(d-1)^{h+1-q} - 1] \sum_{k \in V_j^q} \boldsymbol{\delta}_k = [(d-1)^{h+2-n} - 1]\mathbf{x}_j - [(d-1)^{h+1-n} - 1]\mathbf{x}_{p(j)},$$

and thus (dividing by $d-2$) we obtain

$$\sum_{q=n}^{h} \theta(d, h+1-q) \sum_{k \in V_j^q} \boldsymbol{\delta}_k = \theta(d, h+2-n)\mathbf{x}_j - \theta(d, h+1-n)\mathbf{x}_{p(j)}, \tag{21}$$

which implies

$$\theta(d, h+2-n)\bar{\mathbf{x}}_j - \theta(d, h+1-n)\bar{\mathbf{x}}_{p(j)} = 0.$$

$\square$

**Proposition 7.2** *Let $\mathbf{x}_0$ be the generator corresponding to the root 0. Then $\operatorname{ord}(\bar{\mathbf{x}}_0) = d(d-1)^h$.*

**Proof:** Let $j \in C_0 = \{1, \ldots, d\}$. By equation (21) for $n := 1$ we have:

$$\sum_{q=1}^{h} \theta(d, h+1-q) \sum_{k \in V_j^q} \boldsymbol{\delta}_k = \theta(d, h+1)\mathbf{x}_j - \theta(d, h)\mathbf{x}_0.$$

Adding the above $d$ equations we obtain

$$\sum_{q=1}^{h} \theta(d, h+1-q) \sum_{k \in S_q} \boldsymbol{\delta}_k = \theta(d, h+1) \sum_{j \in C_0} \mathbf{x}_j - d\theta(d, h)\mathbf{x}_0. \qquad (22)$$

From the definition of $\boldsymbol{\delta}_0$ we obtain

$$\theta(d, h+1)\boldsymbol{\delta}_0 = \theta(d, h+1)(d\mathbf{x}_0 - \sum_{j \in C_0} \mathbf{x}_j). \qquad (23)$$

Adding equation (23) to equation (22) and using the definition of the $\theta(d,n)$ (Formula (16)) we obtain

$$\sum_{q=0}^{h} \theta(d, h+1-q) \sum_{k \in S_q} \boldsymbol{\delta}_k = [d/(d-2)][(d-1)^{h+1} - 1 - (d-1)^h + 1]\mathbf{x}_0 \Rightarrow$$

$$\sum_{q=0}^{h} \theta(d, h+1-q) \sum_{k \in S_q} \boldsymbol{\delta}_k = d(d-1)^h \mathbf{x}_0.$$

By Lemma 3.1 (noting that $\theta(d,1) = 1$) we obtain $\operatorname{ord}(\bar{\mathbf{x}}_0) = d(d-1)^h$. $\square$

**Lemma 7.3** *Let $\ell$ be a leaf. Then*

$$\operatorname{ord}(\bar{\mathbf{x}}_\ell) \mid (d-1)^h \cdot \operatorname{lcm}\{d\theta(d, h+1), \theta(d, h), \theta(d, h-1), \ldots, \theta(d, 2)\}.$$

**Proof:** Let $\ell = k_h, k_{h-1}, \ldots, k_0 = 0$ be the path from $\ell$ to the root 0. By Lemma 7.1 we have that

$$\theta(d, h-n+1)\bar{\mathbf{x}}_{k_{n+1}} = \theta(d, h-n)\bar{\mathbf{x}}_{k_n} \quad (0 \le n \le h-1),$$

and by Lemma 7.2 we have that

$$\operatorname{ord}(\bar{\mathbf{x}}_0) = d(d-1)^h.$$

Apply Lemma 3.3 with $t := h$, $\mathbf{z}_n := \bar{\mathbf{x}}_{k_n}$ $(0 \le n \le h)$, $r_n := \theta(d, h+1-n)$ $(0 \le n \le h)$ and $s := d(d-1)^h$ to obtain

$$\operatorname{ord}(\bar{\mathbf{x}}_h) \mid \operatorname{lcm}\{(d-1)^h d\theta(d, h+1), \theta(d, h), \ldots, \theta(d, 2)\}. \qquad (24)$$

The numbers $d$ and $\theta(d,n)$ are all relatively prime to $(d-1)^h$ and the conclusion follows. $\square$

**Lemma 7.4** *Let $\ell$ be a leaf. Then the exponent of $G(d,h)$ is the order of $\bar{\mathbf{x}}_\ell$.*

**Proof:** Let $f = f(d,h) := \text{ord}(\bar{\mathbf{x}}_\ell)$ for all leaves $\ell$. We show that for every vertex $i$, $\text{ord}(\bar{\mathbf{x}}_i) \mid f$, using induction on the distance of $i$ from the boundary. If $i \in S_{h-1}$, then $d\bar{\mathbf{x}}_{\ell_0} = \bar{\mathbf{x}}_i$ for some leaf $\ell_0 \in V_i^h$. So $f\bar{\mathbf{x}}_i = 0 \Rightarrow \text{ord}(\bar{\mathbf{x}}_i) \mid f$. Now, let $i \in S_{h-n}$ and let $m \in C_i$. By the group relations $\bar{\mathbf{x}}_i = d\bar{\mathbf{x}}_m - \sum_{j \in C_m} \bar{\mathbf{x}}_j$ and by the inductive hypothesis $\text{ord}(\bar{\mathbf{x}}_m) \mid f$ and $\text{ord}(\bar{\mathbf{x}}_j) \mid f$, $\forall j \in C_m$.
Therefore, $f\bar{\mathbf{x}}_i = 0 \Rightarrow \text{ord}(\bar{\mathbf{x}}_i) \mid f$. Therefore $f$ is the largest order of an element in $G(d,h)$. □

Lemmas 7.3 and 7.4 yield an upper bound for $\exp(d,h)$:

**Proposition 7.5** *The exponent $\exp(d,h)$ of $G(d,h)$ satisfies*

$$\exp(d,h) \mid (d-1)^h \cdot \text{lcm}\{d\theta(d,h+1), \theta(d,h), \theta(d,h-1), \ldots, \theta(d,2)\}.$$

**Proposition 7.6** *Let $j_1, j_2 \in S_n$ be siblings ($1 \leq n \leq h$). Then the order of $\bar{\mathbf{x}}_{j_1} - \bar{\mathbf{x}}_{j_2}$ is $\theta(d, h+2-n)$.*

**Proof:** Using equation (21) for $j_1$ and $j_2$, we obtain:

$$\theta(d,h+2-n)(\mathbf{x}_{j_1} - \mathbf{x}_{j_2}) = \sum_{q=n}^{h} \theta(d,h+1-q)\left[\sum_{k \in V_{j_1}^q} \boldsymbol{\delta}_k - \sum_{k \in V_{j_2}^q} \boldsymbol{\delta}_k\right].$$

By Lemma 3.1 (noting that $\theta(d,1) = 1$) we obtain $\text{ord}(\bar{\mathbf{x}}_{j_1} - \bar{\mathbf{x}}_{j_2}) = \theta(d, h+2-n)$. □

**Proof of Theorem 2.3:** By Propositions 7.2 and 7.6 we have

$$(d-1)^h \text{lcm}\{\theta(d,n) : n = 1, \ldots, h\} \mid \exp(d,h)$$

Now we need to show that $d\theta(d, h+1) \mid \exp(d,h)$. Let $p \mid d\theta(d, h+1)$.

If $p \nmid d$ then $e_p(d\theta(d,h+1)) = e_p(\theta(d,h+1))$ and the result follows from Proposition 7.6. If $p \nmid \theta(d,h+1)$ then $e_p(d\theta(d,h+1)) = e_p(d)$ and the result follows from Proposition 7.2. If $p \mid d$ and $p \mid \theta(d,h+1)$, in the Sylow $p$-subgroup of $G(d,h)$ consider the equation obtained by Lemma 7.1 for $n := h - 1$:

$$e_p(\theta(d,h+1))\bar{\mathbf{x}}_1 = e_p(\theta(d,h))\bar{\mathbf{x}}_0. \tag{25}$$

We have $p \nmid \theta(d,h)$. By equation (25) and by Lemma 7.2 we obtain $e_p(\bar{\mathbf{x}}_1) = e_p(d\theta(d,h+1))$. □

The next Proposition is interesting in its own right:

**Proposition 7.7** *Let $\mathbf{y}_n = \mathbf{y}(d,n) := \sum_{i \in S_n} \mathbf{x}_i$ ($1 \leq n \leq h$). Then the order of $\bar{\mathbf{y}}_n$ is $(d-1)^{h+1-n}$.*

**Proof:** Adding all the vectors $\boldsymbol{\delta}_k$ yields

$$(d-1)\mathbf{y}_h = \sum_{q=0}^{h} \sum_{k \in S_q} \boldsymbol{\delta}_k.$$

By Lemma 3.1 we have $\operatorname{ord}(\bar{\mathbf{y}}_h) = d-1$.

We proceed by induction on the distance from the boundary. Assume

$$\operatorname{ord}(\bar{\mathbf{y}}_{n+1}) = (d-1)^{h-n} \tag{26}$$

and observe that

$$(d-1)\mathbf{y}_n - \mathbf{y}_{n+1} = \sum_{q=0}^{n} \sum_{k \in S_q} \boldsymbol{\delta}_k,$$

so

$$(d-1)\bar{\mathbf{y}}_n = \bar{\mathbf{y}}_{n+1}. \tag{27}$$

By equations (26) and (27) we obtain $\operatorname{ord}(\bar{\mathbf{y}}_n) = (d-1)^{h-n+1}$. $\square$

## 8 The order and the $(d-1)$ Hall-subgroup of $G(d,h)$

In this section we prove the order formula (2):

$$|G(d,h)| = d(d-1)^h [\theta(d, h+1)]^{d-1} \prod_{n=1}^{h-1} [\theta(d, h+1-n)]^{(d-2)d(d-1)^{n-1}}.$$

Moreover, we prove that the $(d-1)$ Hall-subgroup of $G(d,h)$ is cyclic of order $(d-1)^h$ (Theorem 2.8).

**Definition 8.1** Let $G_n(d,h)$ be the subgroup of $G(d,h)$ generated by $\{\bar{\mathbf{x}}_i : i \in B_n\}$ $(0 \leq n \leq h)$.

**Notation:** Let $B \subset V(\mathcal{T}(d,h))$. Denote by $\operatorname{pr}_B$ the projection map from $\mathbb{Z}^{V(\mathcal{T}(d,h))}$ to $\oplus_{i \in B} \mathbb{Z}\mathbf{x}_i$. So, for $\mathbf{v} = (v_i) \in \mathbb{Z}^{V(\mathcal{T}(d,h))}$ we have $\operatorname{pr}_B(\mathbf{v}) := (v_i)_{i \in B}$.

**Definition 8.2** Let $j \in S_n$ $(0 \leq n \leq h-1)$.

(i) Let $\mathcal{S}_j(d,h)$ denote the subgraph of $\mathcal{T}(d,h)$ induced on the set $V_j$ of descendants of $j$. Note that $\mathcal{S}_0(d,h) = \mathcal{T}(d,h)$. Adjoin a "sink" to $V_j = V(\mathcal{S}_j(d,h))$ and if $n \geq 1$, then join $j$ and every leaf of $\mathcal{S}_j(d,h)$ with 1 and $d-1$ edges respectively to the sink. Let $\hat{\mathcal{S}}_j(d,h)$ denote this augmented graph. If $n = 0$, i.e., if $j$ is the root 0, let $\hat{\mathcal{S}}_j(d,h) := \hat{\mathcal{T}}(d,h)$.

(ii) Let $M_j(d,h)$ be the sandpile group of $\hat{\mathcal{S}}_j(d,h)$.

(iii) Let $K_j(d,h) := M_j(d,h)/\langle \bar{\mathbf{x}}_j \rangle$, where $\langle \bar{\mathbf{x}}_j \rangle$ is the cyclic subgroup of $M_j(d,h)$ generated by $\bar{\mathbf{x}}_j$, where, for the purposes of this definition, $\bar{\mathbf{x}}_j$ denotes the image of $\mathrm{pr}_{V_j}(\mathbf{x}_j)$ in the sandpile group $M_j(d,h)$.

**Lemma 8.3** *Let $j \in S_n$ ($0 \le n \le h$). The vectors $\{\mathrm{pr}_{V_j}(\boldsymbol{\delta}_k) : k \in V_j\}$ are linearly independent over $\mathbb{Q}$.*

**Proof:** The vectors $\{\mathrm{pr}_{V_j}(\boldsymbol{\delta}_k) : k \in V_j\}$ are the row vectors of the reduced Laplacian matrix of $\hat{\mathcal{S}}_j(d,h)$. Therefore

$$M_j(d,h) \simeq \bigoplus_{k \in V_j} \mathbb{Z}\mathrm{pr}_{V_j}(\mathbf{x}_k) \Big/ \sum_{k \in V_j} \mathbb{Z}\mathrm{pr}_{V_j}(\boldsymbol{\delta}_k), \tag{28}$$

and since $M_j(d,h)$ is finite the result follows. □

**Definition 8.4** *Let $j$ be a vertex and let $B \subset V(\mathcal{T}(d,h))$. Let $\Gamma(j, B)$ be the vector subspace of $\mathbb{Q}^B$ generated by the $\{\mathrm{pr}_B(\boldsymbol{\delta}_k) : k \in V_j\}$.*

**Lemma 8.5** *Let $j \in S_n$ ($0 \le n \le h$) and let $A \subset V_j$, $|A| = |V_j| - 1$. Then $\{\mathrm{pr}_{V_j \setminus \{j\}}(\boldsymbol{\delta}_k) : k \in A\}$ is linearly independent over $\mathbb{Q}$.*

**Proof:** We have

$$\Gamma(j, V_j \setminus \{j\}) \simeq \Gamma(j, V_j)/\mathbb{Q}\mathrm{pr}_{V_j}(\mathbf{x}_j) \ (0 \le n \le h).$$

By Lemma (8.3) we have

$$\dim(\Gamma(j, V_j)) = |V(j)| \ (0 \le n \le h), \tag{29}$$

so

$$\dim(\Gamma(j, V_j \setminus \{j\})) = |V_j| - 1. \tag{30}$$

Recall that $C_j$ is the set of children of $j$. Now let us consider equation (21). Let us replace $j$ by $m$ and $n$ by $n+1$; and let us apply $\mathrm{pr}_{V_j \setminus \{j\}}$ to (21) with this modified notation. Noting that $p(m) = j$ and therefore $\mathrm{pr}_{V_j \setminus \{j\}}(\mathbf{x}_{p(m)}) = \mathrm{pr}_{V_j \setminus \{j\}}(\mathbf{x}_j) = 0$ we then obtain

$$\theta(d, h+1-n)\mathrm{pr}_{V_j \setminus \{j\}}(\mathbf{x}_m) = \sum_{q=n+1}^{h} \theta(d, h+1-q)\left[\sum_{k \in V_m^q} \mathrm{pr}_{V_j \setminus \{j\}}(\boldsymbol{\delta}_k)\right] \ (m \in C_j). \tag{31}$$

Also, by the definition of $\delta_j$ we have

$$\mathrm{pr}_{V_j \setminus \{j\}}(\boldsymbol{\delta}_j) + \sum_{m \in C_j} \mathrm{pr}_{V_j \setminus \{j\}}(\mathbf{x}_m) = 0. \tag{32}$$

Multiplying equation (32) by $\theta(d, h + 1 - n)$ and using equation (31) for all $m \in C_j$ we obtain

$$\theta(d,h+1-n)\mathrm{pr}_{V_j\setminus\{j\}}(\boldsymbol{\delta}_j) + \sum_{m\in C_j}\sum_{q=n+1}^{h}\theta(d,h+1-q)\left[\sum_{k\in V_m^q}\mathrm{pr}_{V_j\setminus\{j\}}(\boldsymbol{\delta}_k)\right] = 0. \quad (33)$$

Noting that all the coefficients $\theta(d, h + 1 - q)$ in (33) are non-zero, equations (30) and (33) show that omitting any of the vectors from the set $\{\mathrm{pr}_{V_j\setminus\{j\}}(\boldsymbol{\delta}_k) : k \in V_j\}$, the rest is a basis of $\Gamma(j, V_j\setminus\{j\})$. □

**Definition 8.6** Let $j \in V(\mathcal{T}(d,h))$. Let us pick a leaf $\ell_j \in V_j^h$. Set $U_j := V_j\setminus\{\ell_j\}$.

The next Lemma describes $|V_j^h|$ distinct sets of relations, each of which defines $K_j(d, h)$.

**Lemma 8.7** Let $j \in S_n$ $(0 \leq n \leq h - 1)$. Then

$$K_j(d,h) \simeq \bigoplus_{k\in V_j\setminus\{j\}}\mathbb{Z}\mathrm{pr}_{V_j\setminus\{j\}}(\mathbf{x}_k) \bigg/ \sum_{k\in U_j}\mathbb{Z}\mathrm{pr}_{V_j\setminus\{j\}}(\boldsymbol{\delta}_k).$$

**Proof:** By equation (28) and by the definition of $K_j(d, h)$ we have

$$K_j(d,h) \simeq \bigoplus_{k\in V_j}\mathbb{Z}\mathrm{pr}_{V_j}(\mathbf{x}_k) \bigg/ \left[\sum_{k\in V_j}\mathbb{Z}\mathrm{pr}_{V_j}(\boldsymbol{\delta}_k) + \mathbb{Z}\mathrm{pr}_{V_j}(\mathbf{x}_j)\right] \simeq$$

$$\simeq \bigoplus_{k\in V_j\setminus\{j\}}\mathbb{Z}\mathrm{pr}_{V_j\setminus\{j\}}(\mathbf{x}_k) \bigg/ \sum_{k\in V_j}\mathbb{Z}\mathrm{pr}_{V_j\setminus\{j\}}(\boldsymbol{\delta}_k). \quad (34)$$

By equation (33) (noting that $\theta(d, 1) = 1$) we have

$$\sum_{k\in V_j}\mathbb{Z}\mathrm{pr}_{V_j\setminus\{j\}}(\boldsymbol{\delta}_k) = \sum_{k\in U_j}\mathbb{Z}\mathrm{pr}_{V_j\setminus\{j\}}(\boldsymbol{\delta}_k). \quad (35)$$

Now combine equation (34) with equation (35). □

**Lemma 8.8** $G(d,h)/G_n(d,h) \simeq \bigoplus_{j\in S_n} K_j(d,h)$ $(0 \leq n \leq h - 1)$.

**Proof:** Set $V := V(\mathcal{T}(d, h))$. We have

$$G(d,h)/G_n(d,h) \simeq \bigoplus_{j\in V}\mathbb{Z}\mathbf{x}_j \bigg/ \left[\sum_{j\in V}\mathbb{Z}\boldsymbol{\delta}_j + \sum_{j\in B_n}\mathbb{Z}\mathbf{x}_j\right] \simeq$$

$$\simeq \bigoplus_{j\in V\setminus B_n}\mathbb{Z}\mathrm{pr}_{V\setminus B_n}(\mathbf{x}_j) \bigg/ \left[\sum_{j\in V\setminus B_{n-1}}\mathbb{Z}\mathrm{pr}_{V\setminus B_n}(\boldsymbol{\delta}_j)\right]. \quad (36)$$

Let $j_1, j_2 \in S_n$ with $j_1 \neq j_2$ and let $j_1 \preceq k_1$ and $j_2 \preceq k_2$. Then the vectors $\mathrm{pr}_{V \setminus B_n}(\boldsymbol{\delta}_{k_1})$ and $\mathrm{pr}_{V \setminus B_n}(\boldsymbol{\delta}_{k_2})$ are disjoint. Therefore by equation (36) we obtain

$$G(d,h)/G_n(d,h) \simeq \bigoplus_{j \in S_n} \left[ \bigoplus_{k \in V_j \setminus \{j\}} \mathbb{Z}\mathrm{pr}_{V_j \setminus \{j\}}(\mathbf{x}_k) \Big/ \sum_{k \in V_j \setminus \{j\}} \mathbb{Z}\mathrm{pr}_{V_j \setminus \{j\}}(\boldsymbol{\delta}_k) \right]. \tag{37}$$

By equations (34) and (37) we obtain the result. $\square$

**Lemma 8.9** *Let $j \in S_n$ and let $i := p(j)$ be the parent of $j$. Then the order of $\mathrm{pr}_{V_i \setminus \{i\}}(\mathbf{x}_j)$ in $K_i(d,h)$ is $\theta(d, h+2-n)$ $(1 \leq n \leq h)$.*

**Proof:** Let us apply $\mathrm{pr}_{V_i \setminus \{i\}}$ to equation (21). Noting that $\mathrm{pr}_{V_i \setminus \{i\}}(\mathbf{x}_{p(j)}) = \mathrm{pr}_{V_i \setminus \{i\}}(\mathbf{x}_i) = 0$ we obtain

$$\theta(d, h+2-n)\mathrm{pr}_{V_i \setminus \{i\}}(\mathbf{x}_j) = \sum_{q=n}^{h} \theta(d, h+1-q) \left[ \sum_{k \in V_j^q} \mathrm{pr}_{V_i \setminus \{i\}}(\boldsymbol{\delta}_k) \right]. \tag{38}$$

In equation (38) we have obtained an expression of $\theta(d, h+2-n)\mathrm{pr}_{V_i \setminus \{i\}}(\mathbf{x}_j)$ as an integer linear combination of the vectors $\{\mathrm{pr}_{V_i \setminus \{i\}}(\boldsymbol{\delta}_k) : k \in V_j\}$. Let $j'$ be a sibling of $j$ and let $\ell_{j'} \in V_{j'}^h \subset V_i^h$ be a leaf below $j'$. Then $V_j \subset V_i \setminus \{\ell_{j'}\}$. By Lemma 8.7 we have that $\{\mathrm{pr}_{V_i \setminus \{i\}}(\boldsymbol{\delta}_k) : k \in V_j\}$ is part of a lattice basis that defines $K_i(d,h)$. Noting that $\theta(d,1) = 1$, the conclusion follows from Lemma 3.1. $\square$

**Proposition 8.10** *Let $j \in S_n$. Then the order of $\mathbf{x}_j$ in $G(d,h)/G_{n-1}(d,h)$ is $\theta(d, h+2-n)$ $(1 \leq n \leq h)$.*

**Proof:** The conclusion follows from Lemmas 8.8 and 8.9. $\square$

**Lemma 8.11** *Let $j \in S_n$ $(0 \leq n \leq h-1)$. Let $C_j' \subset C_j$ with $|C_j'| = |C_j| - 1$. Assume*

$$\sum_{m \in C_j'} r_m \mathrm{pr}_{V_j \setminus \{j\}}(\mathbf{x}_m) = 0 \text{ in } K_j(d,h). \tag{39}$$

*Then $r_m \mathrm{pr}_{V_j \setminus \{j\}}(\mathbf{x}_m) = 0$ in $K_j(d,h)$ $(m \in C_j')$.*

**Proof:** Assume for a contradiction that for some $m_0 \in C_j'$ we have

$$r_{m_0} \mathrm{pr}_{V_j \setminus \{j\}}(\mathbf{x}_{m_0}) \neq 0 \text{ in } K_j(d,h). \tag{40}$$

Lemma 8.9 applied with $j$ replaced by $m_0$, $i$ replaced by $j$ and $n$ replaced by $n+1$ implies that

$$\theta(d, h+1-n)\mathrm{pr}_{V_j\setminus\{j\}}(\mathbf{x}_{m_0}) = 0 \text{ in } K_j(d,h). \tag{41}$$

By equations (40) and (41) we have

$$\theta(d, h+1-n) \nmid r_{m_0}. \tag{42}$$

Let $C_j'' := C_j'\setminus\{m_0\}$. By the Euclidean algorithm for $r_{m_0}$ and $\theta(d, h+1-n)$ and by equations (39) and (41) we obtain

$$\sum_{m \in C_j'} s_m \mathrm{pr}_{V_j\setminus\{j\}}(\mathbf{x}_m) = 0 \text{ in } K_j(d,h), \tag{43}$$

where

$$s_{m_0} = \gcd(r_{m_0}, \theta(d, h+1-n)) \text{ and } s_m \in \mathbb{Z}\ (m \in C_j''). \tag{44}$$

Set

$$t_1 := \gcd\{s_m : m \in C_j'\}, \tag{45}$$

and

$$t_2 := \gcd\{s_m\theta(d, h+1-q) : m \in C_j', q = n+1, \ldots, h\}. \tag{46}$$

**Claim 8.12** $t_1 = t_2$.

**Proof of Claim 8.12:** Clearly $t_1 \mid t_2$. For the other direction consider the definition (46). Setting $q := h$ and recalling that $\theta(d, 1) = 1$ we obtain

$$t_2 \mid s_m\ (m \in C_j'),$$

and therefore $t_2 \mid t_1$. □

Observe that $s_{m_0} \mid \theta(d, h+1-n)$ and therefore $t_1 \mid \theta(d, h+1-n)$.

**Claim 8.13** *The order of $\sum_{m \in C_j'} s_m \mathrm{pr}_{V_j\setminus\{j\}}(\mathbf{x}_m)$ in $K_j(d,h)$ is equal to $\theta(d, h+1-n)/t_1$.*

**Proof of Claim 8.13:** We have

$$\theta(d, h+1-n)/t_1 \sum_{m \in C_j'} s_m \mathrm{pr}_{V_j\setminus\{j\}}(\mathbf{x}_m) =$$

$$= \sum_{m \in C_j'} (s_m/t_1)\theta(d, h+1-n)\mathrm{pr}_{V_j\setminus\{j\}}(\mathbf{x}_m). \tag{47}$$

Substituting $\theta(d, h+1-n)\text{pr}_{V_j\setminus\{j\}}(\mathbf{x}_m)$ from equation (31) into (47) we obtain

$$\theta(d, h+1-n)/t_1 \sum_{m \in C_j'} s_m \text{pr}_{V_j\setminus\{j\}}(\mathbf{x}_m) =$$

$$= \sum_{m \in C_j'} s_m/t_1 \sum_{q=n+1}^{h} \theta(d, h+1-q) \sum_{k \in V_m^q} \text{pr}_{V_j\setminus\{j\}}(\boldsymbol{\delta}_k) =$$

$$= \sum_{m \in C_j'} \sum_{q=n+1}^{h} s_m \theta(d, h+1-q)/t_1 \sum_{k \in V_m^q} \text{pr}_{V_j\setminus\{j\}}(\boldsymbol{\delta}_k). \quad (48)$$

Let $\{m_1\} = C_j \setminus C_j'$ and let $\ell \in V_{m_1}^h$ be a leaf below $m_1$. In equation (48) we have obtained an expression of the vector $\theta(d, h+1-n)/t_1 \sum_{m \in C_j'} s_m \text{pr}_{V_j\setminus\{j\}}(\mathbf{x}_m)$ as a linear combination of the vectors

$$\{\text{pr}_{V_j\setminus\{j\}}(\boldsymbol{\delta}_k) : k \in \cup_{m \in C_j'} V_m\} \subset \{\text{pr}_{V_j\setminus\{j\}}(\boldsymbol{\delta}_k) : k \in V_j \setminus \{\ell\}\}.$$

By Lemma 8.7 the set $\{\text{pr}_{V_j\setminus\{j\}}(\boldsymbol{\delta}_k) : k \in \cup_{m \in C_j'} V_m\}$ is part of a lattice basis that defines $K_j(d, h)$. Also, by Claim 8.12 we have

$$\gcd\{s_m \theta(d, h+1-q)/t_1 : m \in C'j, q = n+1, \ldots, h\} =$$
$$= (1/t_1)\gcd\{s_m \theta(d, h+1-q) : m \in C'j, q = n+1, \ldots, h\} = t_2/t_1 = 1,$$

so, the coefficients in equation (48) are relatively prime. The conclusion follows by Lemma 3.1. □

By equation (43) and by Claim 8.13 we conclude

$$\theta(d, h+1-n) = t_1 = \gcd\{s_m : m \in C_j'\}.$$

Therefore $\theta(d, h+1-n) \mid s_{m_0}$ and hence by equation (44), $\theta(d, h+1-n) \mid r_{m_0}$, which contradicts our original assumption (equation (42)). □

**Definition 8.14** Let $j \in S_n$ ($0 \le n \le h-1$). Let $E_j(d, h)$ be the subgroup of $K_j(d, h)$ generated by the images of the $\{\text{pr}_{V_j\setminus\{j\}}(\mathbf{x}_m) : m \in C_j\}$ in $K_j(d, h)$ represented as in Lemma 8.7.

**Lemma 8.15** Let $j \in S_n$. Then

$$E_j \simeq \underbrace{\mathbb{Z}_{\theta(d,h+1-n)} \oplus \cdots \oplus \mathbb{Z}_{\theta(d,h+1-n)}}_{|C_j|-1 \text{ times}} \quad (0 \le n \le h-1).$$

**Proof:** By Lemma 8.7 we have that $\text{pr}_{V_j\setminus\{j\}}(\boldsymbol{\delta}_j)$ is a member of a lattice basis that defines $K_j(d, h)$. Therefore

$$\text{pr}_{V_j\setminus\{j\}}(\boldsymbol{\delta}_j) = 0 \text{ in } K_j(d, h). \quad (49)$$

By equations (32) and (49) we obtain

$$\sum_{m \in C_j} \mathrm{pr}_{V_j \setminus \{j\}}(\mathbf{x}_m) = 0 \text{ in } K_j(d,h). \tag{50}$$

Equation (50) combined with Lemmas 8.9 (applied with $j$ replaced by $m$, $i$ replaced by $j$ and $n$ replaced by $n+1$) and 8.11 yield the result. $\square$

**Lemma 8.16**

$$G_{n+1}(d,h)/G_n(d,h) \simeq \bigoplus_{j \in S_n} E_j \ (0 \leq n \leq h-1).$$

**Proof:** The result follows from Lemma 8.8. $\square$

**Lemma 8.17**

$$|G_{n+1}(d,h)/G_n(d,h)| =$$
$$= \begin{cases} \theta(d, h+1)^{d-1} & \text{if } n = 0, \\ \theta(d, h+1-n)^{(d-2)d(d-1)^{n-1}} & \text{if } 1 \leq n \leq h-1. \end{cases}$$

**Proof:** The conclusion follows from Lemmas 8.15 and 8.16. $\square$

**Proof of Theorem 2.5**: Consider the chain of subgroups $G_0(d,h) \leq G_1(d,h) \leq \ldots \leq G_h(d,h) = G(d,h)$. Clearly,

$$|G(d,h)| = |G_0(d,h)| \prod_{n=0}^{h-1} |G_{n+1}(d,h)/G_n(d,h)|.$$

Proposition 7.2 and Lemma 8.17 give the result. $\square$

**Lemma 8.18** *$G(d,h)$ has a cyclic subgroup of order $(d-1)^h$.*

**Proof:** By Proposition 7.2 we have that $\mathrm{ord}(d\bar{\mathbf{x}}_0) = (d-1)^h$. $\square$

**Proof of Theorem 2.8:** We have that $d$ and the numbers $\theta(d,j)$ are relatively prime to $d-1$. This fact combined with Lemma 8.18 and with the order formula (2) (Theorem 2.5) proves Theorem 2.8. $\square$

# 9 Asymptotic evaluation

**Proof of Corollary 2.6:** Formula (2) yields

$$\log_{d-1} g(d,h) = \log_{d-1} d + h + (d-1)\log_{d-1} \theta(d, h+1) +$$
$$+ \sum_{n=1}^{h-1} (d-2)d(d-1)^{n-1} \log_{d-1}\left(\frac{(d-1)^{h+1-n} - 1}{d-2}\right) =$$
$$= \log_{d-1} d + h + (d-1)\log_{d-1} \theta(d, h+1) +$$
$$+ \sum_{n=0}^{h-2} (d-2)d(d-1)^{h-n-2} \log_{d-1}\left(\frac{(d-1)^{n+2} - 1}{d-2}\right).$$

Clearly,

$$\lim_{h\to\infty} \frac{\log_{d-1} g(d,h)}{(d-1)^h} = c_d.$$

□

**Notation:** Let $r, s$ be positive integers. Set $\eta(r,s) := \text{lcm}\{r^n - 1 : n = 1, \ldots, s\}$.

We will need the following estimate (see [12]):

**Theorem 9.1** *For every fixed $r \geq 2$, the following asymptotic equality holds as $s \to \infty$:*

$$\log_r \eta(r,s) \sim \frac{3s^2}{\pi^2}.$$

**Proof of Corollary 2.4:** By formula (1) we obtain

$$\eta(d, h+1) \leq \exp(d, h) \leq d(d-1)^h \eta(d, h+1). \tag{51}$$

The result follows by combining formula (51) with Theorem 9.1 used with parameters $r := d-1$ and $s := h$. □

# 10 Questions for further study

In this paper we have calculated the exponent and the rank of the abelian group $G(d, h)$. It seems plausible that all invariants could be found explicitly, i.e., the group structure can be fully determined. This would be equivalent to knowing the structure of all Sylow subgroups. We now describe an explicit conjecture regarding the rank of each Sylow subgroup.

**Notation:** Let $p$ be a prime. We denote by $S_p(d,h)$ the Sylow $p$-subgroup of $G(d, h)$.

If $p \mid d$ then by Lemma 6.6 we have $\mathrm{rank}(S_p(d,h)) = (d-1)^h$. If $p \mid d-1$ then by Theorem 2.8 we have $\mathrm{rank}(S_p(d,h)) = 1$.

**Notation:**

(i) Let $r, s$ be relatively prime integers with $s > 0$. We denote by $\mathrm{ord}_s(r)$ the multiplicative order of $r \mod s$.

(ii) Let $p$ be a prime with $p \nmid d-1$. We denote by $t_p(d)$ the least positive integer $n$ such that $p \mid \theta(d, n)$. Note that

$$t_p(d) = \begin{cases} p & \text{if } d \equiv 2 \mod p, \\ \mathrm{ord}_p(d-1) & \text{otherwise.} \end{cases}$$

**Conjecture 10.1** *Let $p \nmid d(d-1)$. The rank of the Sylow $p$-subgroup $S_p(d,h)$ of $G(d,h)$ is given by the following formula:*

$$\mathrm{rank}(S_p(d,h)) = \begin{cases} d(d-1)^{h-rt_p(d)} \sum_{q=0}^{r-1} (d-1)^{qt_p(d)} \\ \qquad \text{if } rt_p(d) \leq h \leq (r+1)t_p(d) - 2, \\ d(d-1)^{t_p(d)-1} \sum_{q=0}^{r-1} (d-1)^{qt_p(d)} + d - 1 \\ \qquad \text{if } h = (r+1)t_p(d) - 1. \end{cases} \quad (52)$$